\documentclass[a4paper]{article}
\pdfoutput=1

\usepackage{amsmath, amssymb, amsthm}
\usepackage{mathtools}
\usepackage{mathrsfs}
\usepackage{enumitem}
\usepackage{float}

\newcommand{\F}{\ensuremath{\mathbb{F}}}

\newtheorem{theorem}{Theorem}[section]

\newtheorem{conjecture}[theorem]{Conjecture}
\theoremstyle{definition}

\theoremstyle{remark}
\newtheorem{remark}[theorem]{Remark}

\DeclareMathOperator{\T}{T}
\DeclareMathOperator{\ST}{ST}
\DeclareMathOperator{\GF}{GF}

\DeclarePairedDelimiterX{\Set}[1]{\{}{\}}{
	\,#1\,
}

\title{On the Number of Periodic Points of Quadratic Dynamical Systems Modulo a Prime}
\author{Jakob Streipel}
\date{\today}

\begin{document}

	\maketitle
	
	\begin{abstract}
		
		\noindent
		In 2004 Vasiga and Shallit studied the number of periodic points of two particular discrete quadratic maps modulo prime numbers. They found the asymptotic behaviour of the sum of the number of periodic points for all primes less than some bound, assuming the Extended Riemann Hypothesis. Later that same year Chou and Shparlinski proved this asymptotic result without assuming any unproven hypotheses. 
		
		Inspired by this we perform experiments and find a striking pattern in the behaviour of the sum of the number of periodic points for quadratic maps other than the two particular ones studied previously. From simulations it appears that the sum of the number of periodic points of all quadratic maps of this type behave the same. Finally we find that numerically the distribution of the amounts of periodic points seems to be Rayleigh. 
	
	\end{abstract}

	\section{Introduction} 
	
	Quadratic maps of the type $f_c(x) = x^2 + c$ over finite fields (or rings) are of great interest in both number theory and cryptography. They are used in primality tests such as Lucas--Lehmer and Miller--Rabin, integer factorisation methods like Pollard's $\rho$ method \cite{Pollard1975}, and pseudorandom number generation \cite[Section~4]{BlumBlumShub86}.
	
	In studying these maps we will use the natural definitions and notations. For a map $f : S \to S$ we let $f^r$ define the \emph{$r$-fold composition of $f$} for nonnegative integers $r$, with $f^0$ defined as the identity mapping. We call $f^n$ the $n$th \emph{iterate of $f$}. We will call $x$ a \emph{periodic point (of $f$)} if there exists some positive integer $n$ such that $f^n(x) = x$, with the smallest such $n$ being called the \emph{period of $x$}. 
	
	Moreover we denote by $G_f = (\mathscr{V}, \mathscr{E})$ the \emph{directed graph} of the dynamical system of $f : S \to S$, meaning that the vertices $\mathscr{V}$ is the set of elements of $S$ and its directed edges $\mathscr{E}$ are $(x, f(x))$ for every $x \in S$. By the \emph{orbit of $x$} we mean the directed path in $G_f$  starting at $x$, and with the \emph{tail of $x$} we mean the list of elements $\Set{x, f(x), f^2(x), \ldots}$ before we encounter a periodic point, and the orbit of this periodic point we encounter we call the \emph{cycle of $x$}.
	
	In particular the dynamical system given by $f_0(x) = x^2$ is well-studied and well-known. It has been of great interest ever since Blum et al. introduced the BBS pseudorandom number generator, which uses this mapping modulo the product of two prime numbers. Following this Rogers in \cite{Rogers1996} gives a formula to decompose the graph $G_{f_0}$ into its cyclic components and the trees attached to these cycles, and Hern\'andez et al. in \cite{Hernandez1994} study and completely characterise the orbits of $f_0$.
	
	Similarly well-studied is the map $f_{-2}(x) = x^2 - 2$, which is the one used in the Lucas--Lehmer primality test for Mersenne primes. A BBS-like pseudorandom number generator using this map instead of the original $f_0(x) = x^2$ has been studied by Dur\'an D\'ias and Peinado Dom\'inguez in \cite{Duran2002}. Results regarding the dynamics of this map were obtained by Gilbert et al. in \cite[Section~5]{Gilbert2001}.
	
	Of note is Pollard's $\rho$ method of integer factorisation, making use of any $f_c$, which Pollard in \cite[page~333]{Pollard1975} cautions should \emph{not} be used with $f_{-2}$, as further motivated by \cite[Section~4]{Vasiga2004}.  
	
	Of special interest to us is the results of Vasiga and Shallit in \cite{Vasiga2004}, where in Section 2 and Section 3 they study in great detail the dynamics of $f_0(x) = x^2$ and $f_{-2}(x) = x^2 - 2$, respectively. In doing so they obtain explicit expressions for what we will denote $\T_0(p)$, the number of periodic points of $f_0$ in $\F_p$, the field of $p$ elements, $p$ being prime, as well as the corresponding $\T_{-2}(p)$, being the number of periodic points of $f_{-2}$ in $\F_p$. 
	
	Using these results they go on to study the asymptotic behaviour of
	\[
		\ST_0(N) = \sum_{p \leq N} \T_0(p) \qquad \text{and} \qquad \ST_{-2}(N) = \sum_{p \leq N} \T_{-2}(p),
	\]
	being the sum of the number of periodic points for all primes $p$ less than or equal to some $N$. They find that, assuming the Extended Riemann Hypothesis, the asymptotic behaviour of both of these is
	\[
		\frac{N^2}{6 \log N},
	\]
	which Chou and Shparlinski in \cite{Shparlinski2004} show to be true even without assuming the Extended Riemann Hypothesis. For a detailed working out of Vasiga and Shallit's results on a fairly elementary level, see \cite[Sections~2~\&~3]{Streipel2015}. 
	
	Unfortunately, not as much is known about the number of periodic points $\T_c(p)$ of $f_c(x) = x^2 + c$ over $\F_p$ when $c$ is different from $0$ and $-2$, and it is difficult to apply the methods used for $c = 0$ and $c = -2$---which boil down to expressing the $n$th iterate of $f_c$ as a sum of terms of the form $\eta^{2^n}$, where $\eta$ are some constants---in order to study other $c$.
	
	Having said that, there are certain results pertaining to $\T_c(p)$. Peinado et al. in \cite{Peinado2001} find several different upper bounds on the cycle lengths of $x \mapsto x^2 + c$ over $\F_q$, depending on the nature of the prime power $q$ and the coefficient $c$. They also in \cite[Propositions~3~\&~4]{Peinado2001} find conditions for when this map has cycles of lengths 1 and 2, depending on whether certain expressions in $c$ are quadratic residues in the given field.
	
	There is also \cite[Section~3]{Nilsson2013} in which Nilsson studies the periodic points of all $f_c$ in $\F_p$ for fixed primes $p > 2$. He shows that there are exactly $(p - 1) / 2$ diagonal lines in the so-called Periodic Point Diagram---a means of visualising which combinations of $c$ and $x$ are periodic---that contain no periodic points at all. 
	
	Our approach is to attack the unpredictable behaviour of $\T_c(p)$, $c \neq 0, -2$ from a different angle: instead of attempting to study $\T_c(p)$ directly, we are inspired by the asymptotic results of Vasiga \& Shallit and Chou \& Shparlinski and study experimentally
	\[
		\ST_c(N) = \sum_{p \leq N} \T_c(p),
	\]
	the sum of the number of periodic points of $f_c$ for all primes $p$ less than or equal to some $N$. 
	
	We find a striking pattern in these for $c \neq 0, -2$, and in particular this same pattern appears for all $c \neq 0, -2$ studied. Based on this empirical data we formulate a conjecture describing this behaviour. Finally we study the distribution of $\T_c(p)$, finding that after normalising to $\T_c(p) / \sqrt{p}$ they appear Rayleigh distributed.

	\section{Empirical investigation of $\ST_c(N)$}
	
\begin{figure}
		\centering
		\includegraphics[width = 0.48 \textwidth]{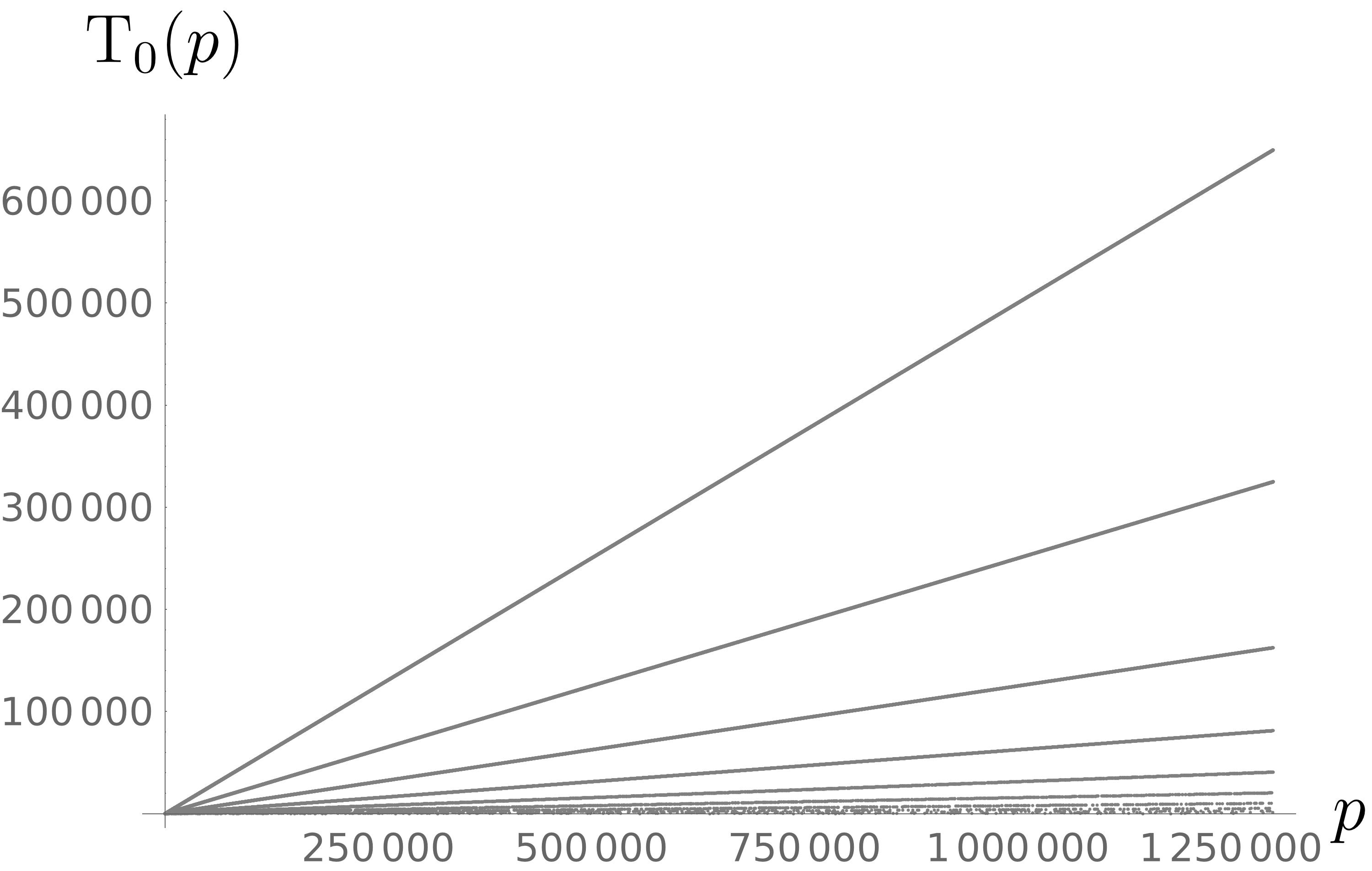} \quad
		\includegraphics[width = 0.48 \textwidth]{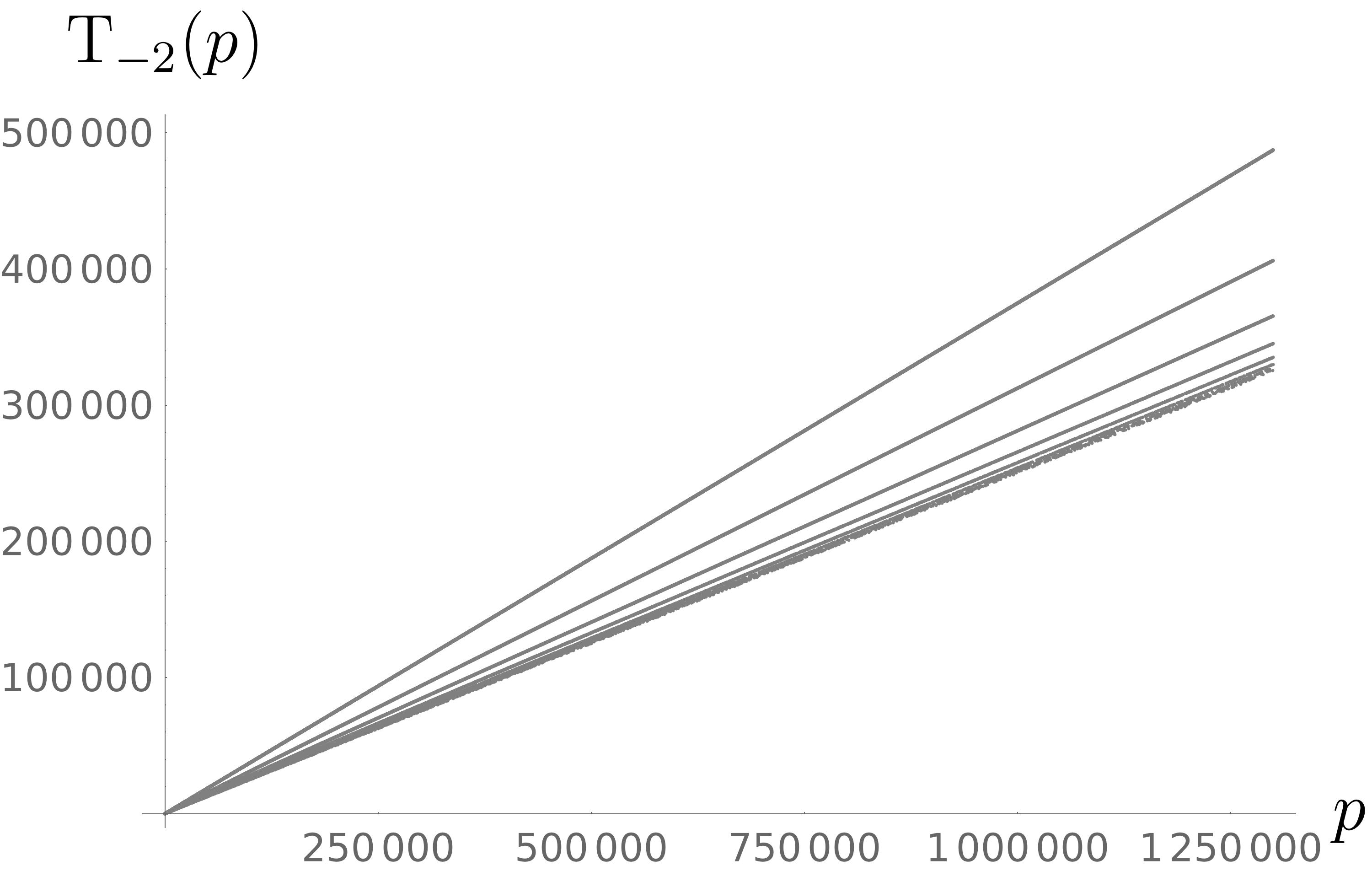}
		
		\includegraphics[width = 0.48 \textwidth]{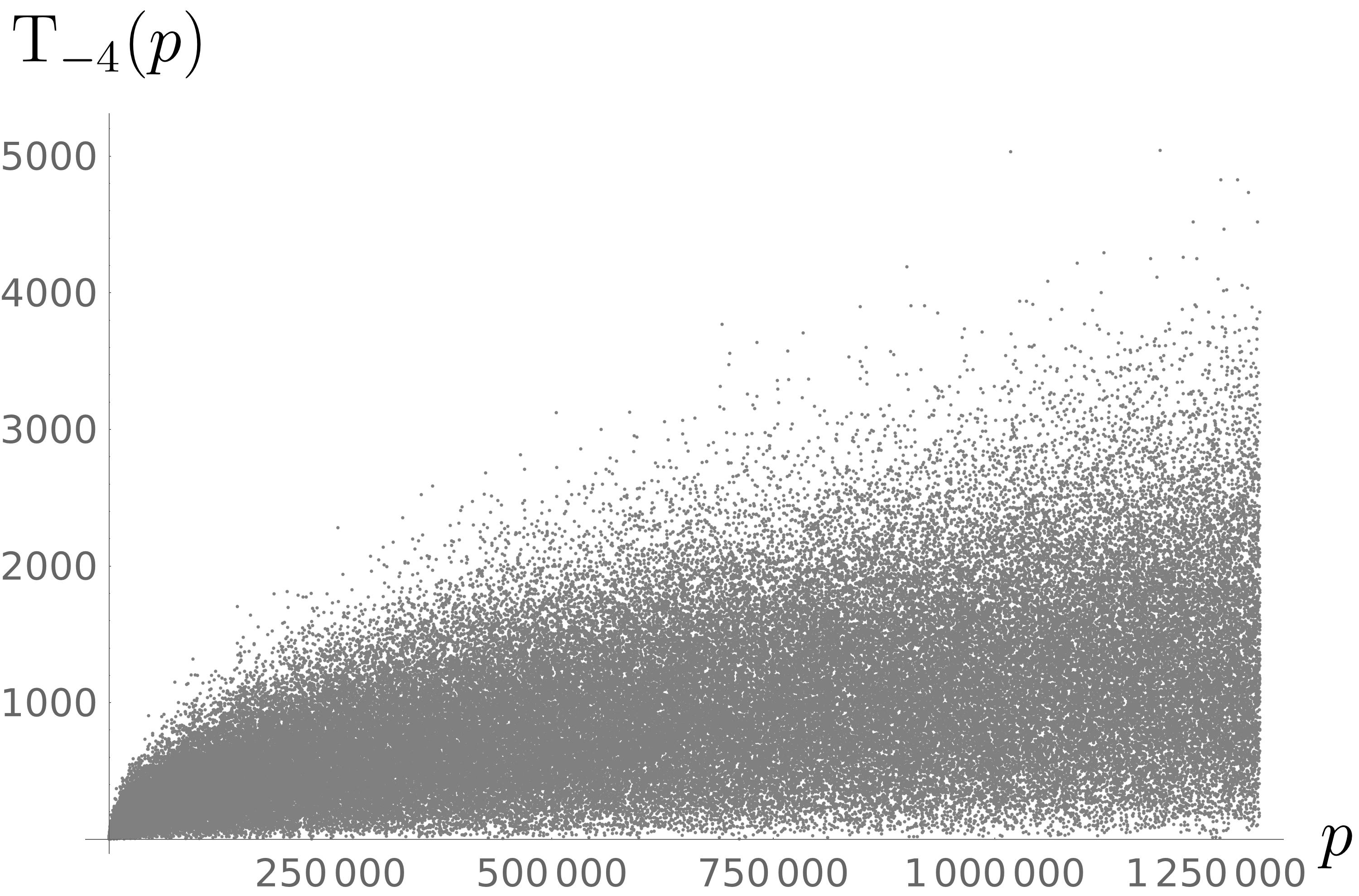} \quad
		\includegraphics[width = 0.48 \textwidth]{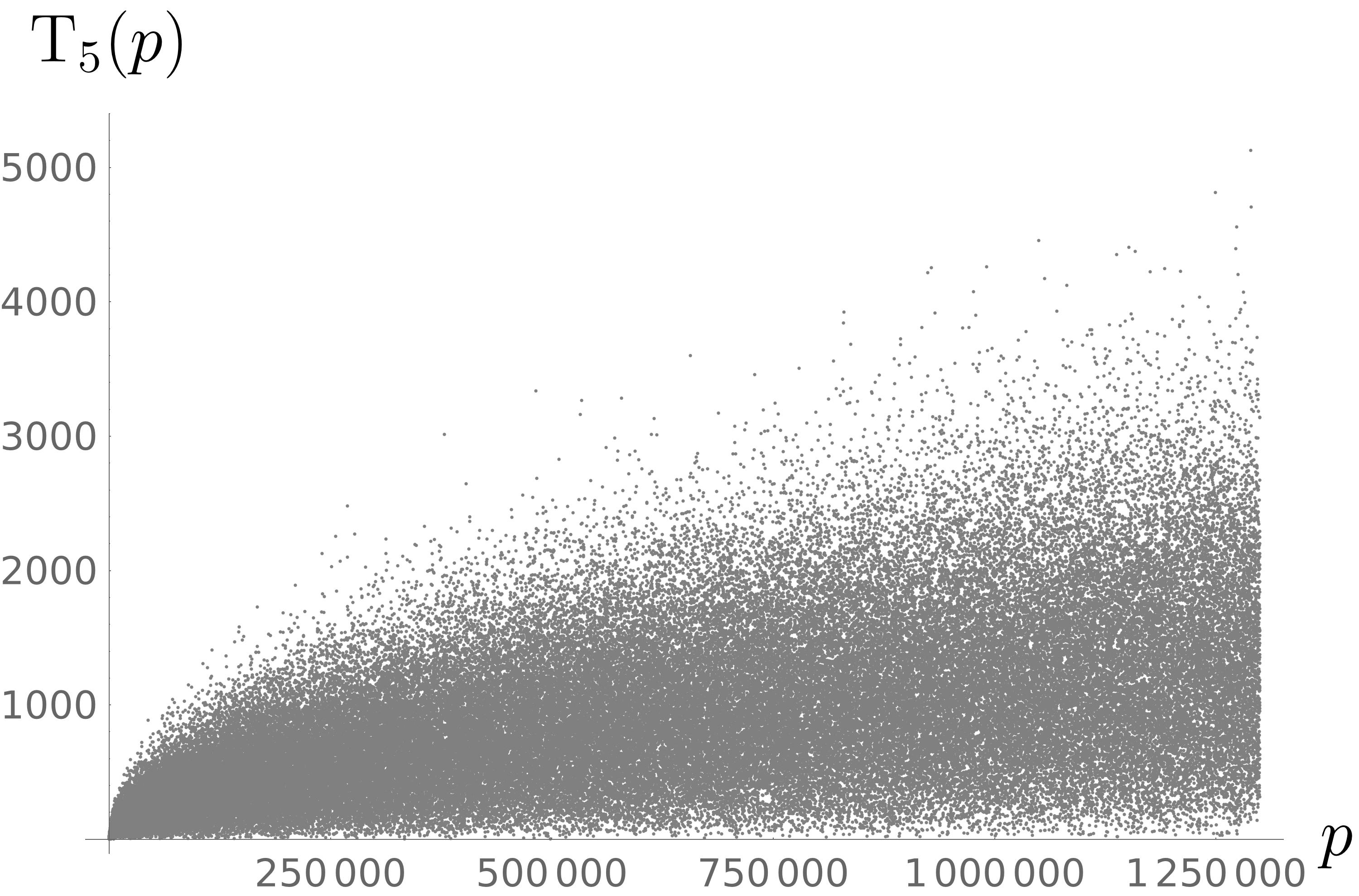}
		\caption{Plots of $\T_c(p)$  for $c = 0, -2, -4, 5$, for the first one hundred thousand prime numbers $p$.}
		\label{fig:Tcpchaos}
	\end{figure}

	We demonstrate the unpredictable behaviour of $\T_c(p)$ for $c \neq 0, -2$ in Figure~\ref{fig:Tcpchaos} by comparing plots of them with plots of $\T_0(p)$ and $\T_{-2}(p)$. It is plain to see that $\T_0(p)$ and $\T_{-2}$ are much more well behaved. Note also that the magnitudes of $\T_0(p)$ and $\T_{-2}(p)$ are far, far greater than that for the other $c$. The plots of other $\T_c(p)$, $c \neq 0, -2$, look similar to that of $\T_{-4}(p)$ and $\T_5(p)$, at least for $c \neq 0, -2$ between $-100$ and $100$.
	
	We compute $\T_c(p)$ using the following fairly na\"ive algorithm:
	\begin{enumerate}
		\item Create a list of all elements in $\F_p$,
		\item Apply $f_c$ on each element in the list,
		\item Remove duplicates from the list,
		\item Repeat from step 2 until the number of elements in the list remains constant. The number of elements in the list is now $\T_c(p)$. 
	\end{enumerate}
	
	\noindent
	The motivation for the algorithm is rather straight forward: if, at any point, two elements in the list are equal, there is no sense in continuing to iterate on both; since they are equal, their orbits are also equal. Moreover, for the number of elements in the list to remain constant, we must have that no two of the remaining elements map to the same element, because then the number of elements in the list would decrease. Thus since all remaining elements must map to different elements, we must be stuck in cycles. 
	
	It is worth noting that this is by no means a particularly fast algorithm, though it is fast enough for our purposes. Computing $\T_c(p)$ for the first ten thousand primes for a fixed $c$ takes in the order of a quarter of an hour on our ordinary personal computer. In addition, since the algorithm needs to store a list of all elements in $\F_p$ it also requires a considerable amount of memory for large primes $p$. Both running time and memory requirements could be improved upon by implementing a tortoise and hare algorithm like Floyd's cycle finding algorithm \cite{Floyd1967} or Brent's algorithm \cite{Brent1980}.
	
	\begin{remark}
		Note that if $c \geq p$, $c$ will be reduced modulo $p$. Therefore, for example $\T_3(5) = \T_{-2}(5)$ since $3 \equiv -2 \pmod{5}$. Because of this there is a certain amount of double counting when considering $\ST_c(N)$, since for example $\ST_3(N)$ and $\ST_{-2}(N)$ for $N \geq 5$ will include the same term $\T_3(5)$. As $N$ increases, this will of course become less and less significant.
	\end{remark}

	\begin{figure}[t!]
		\centering
		\includegraphics[width = 0.48 \textwidth]{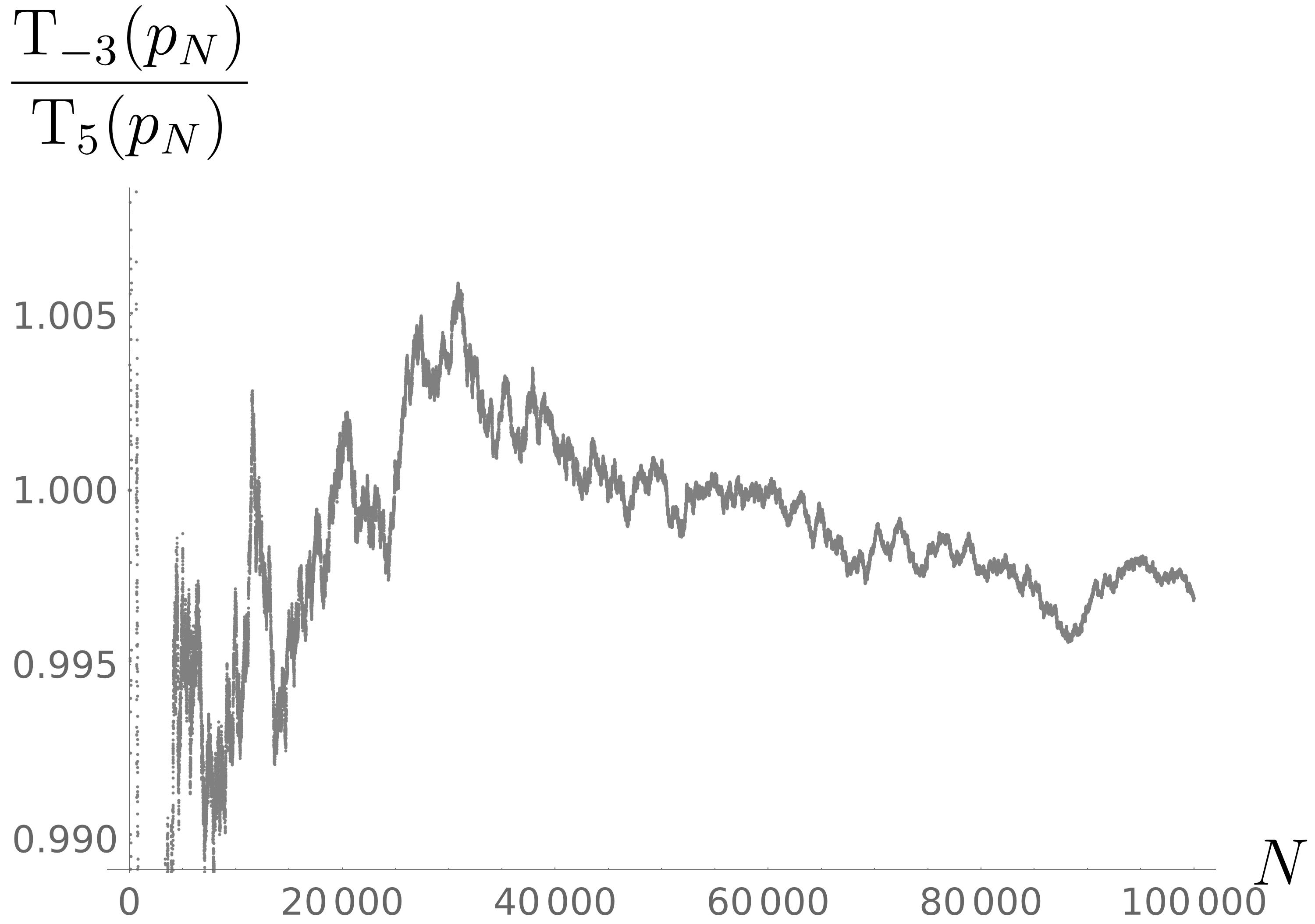} \quad \includegraphics[width = 0.48 \textwidth]{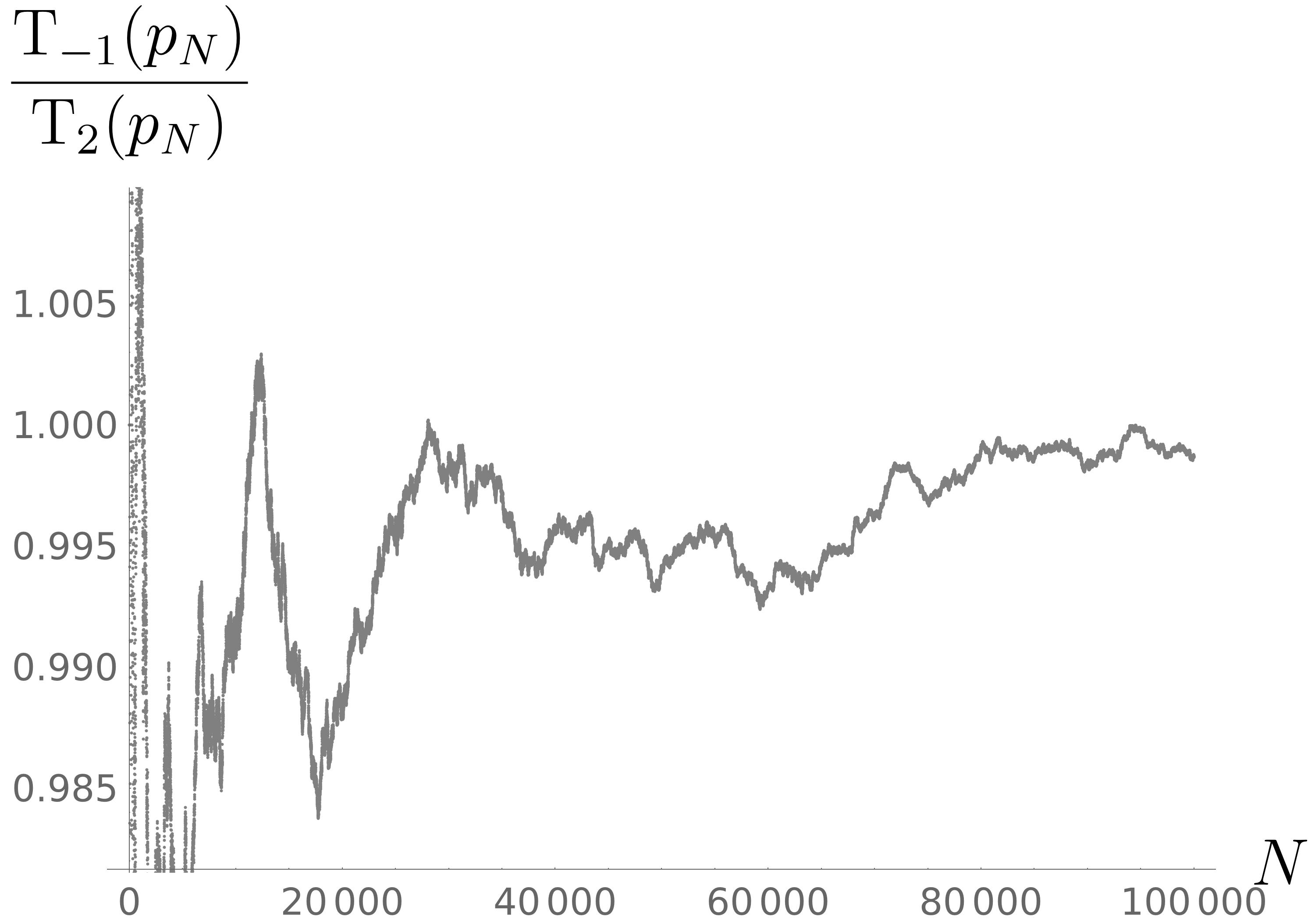}
		
		\includegraphics[width = 0.48 \textwidth]{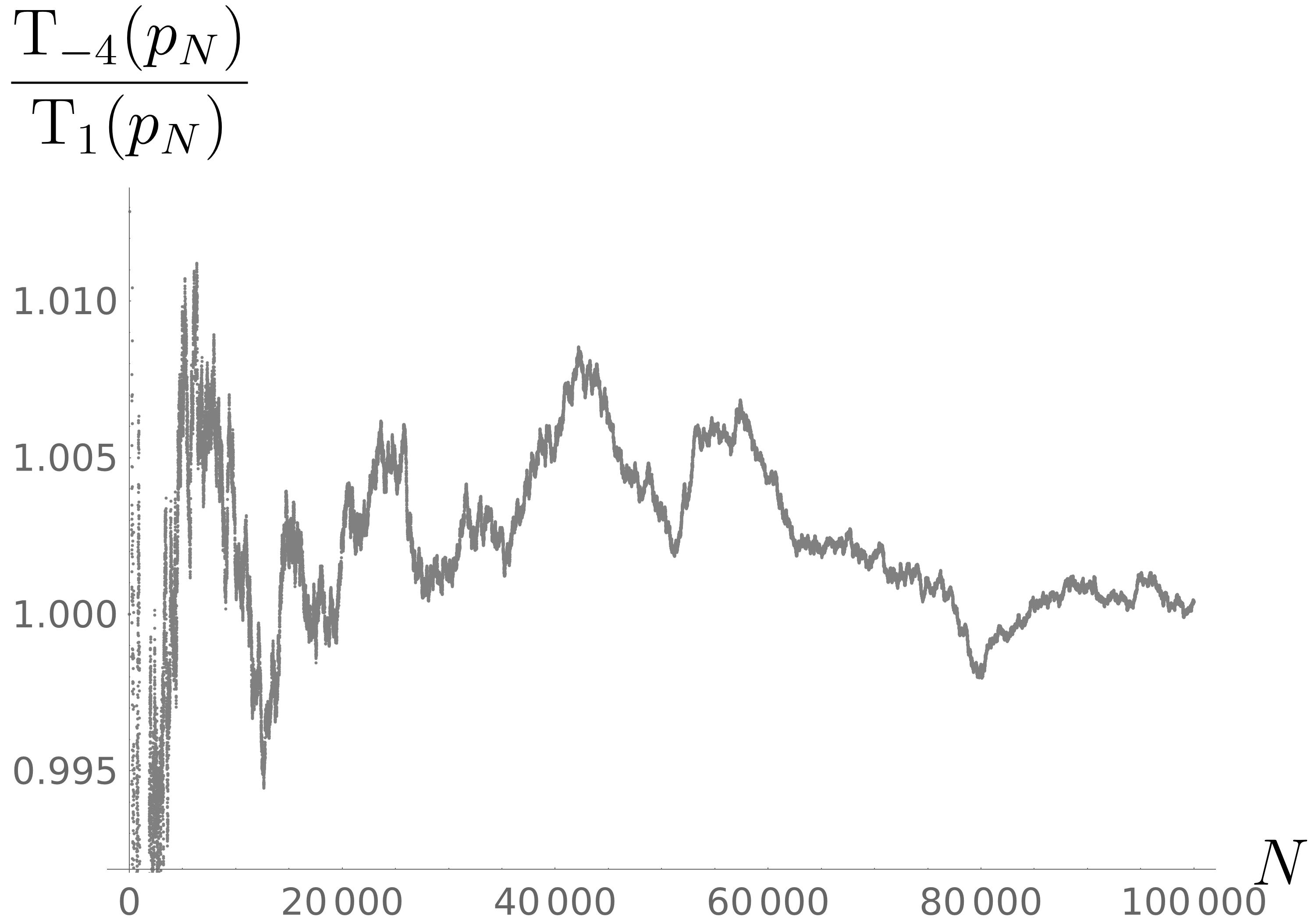} \quad \includegraphics[width = 0.48 \textwidth]{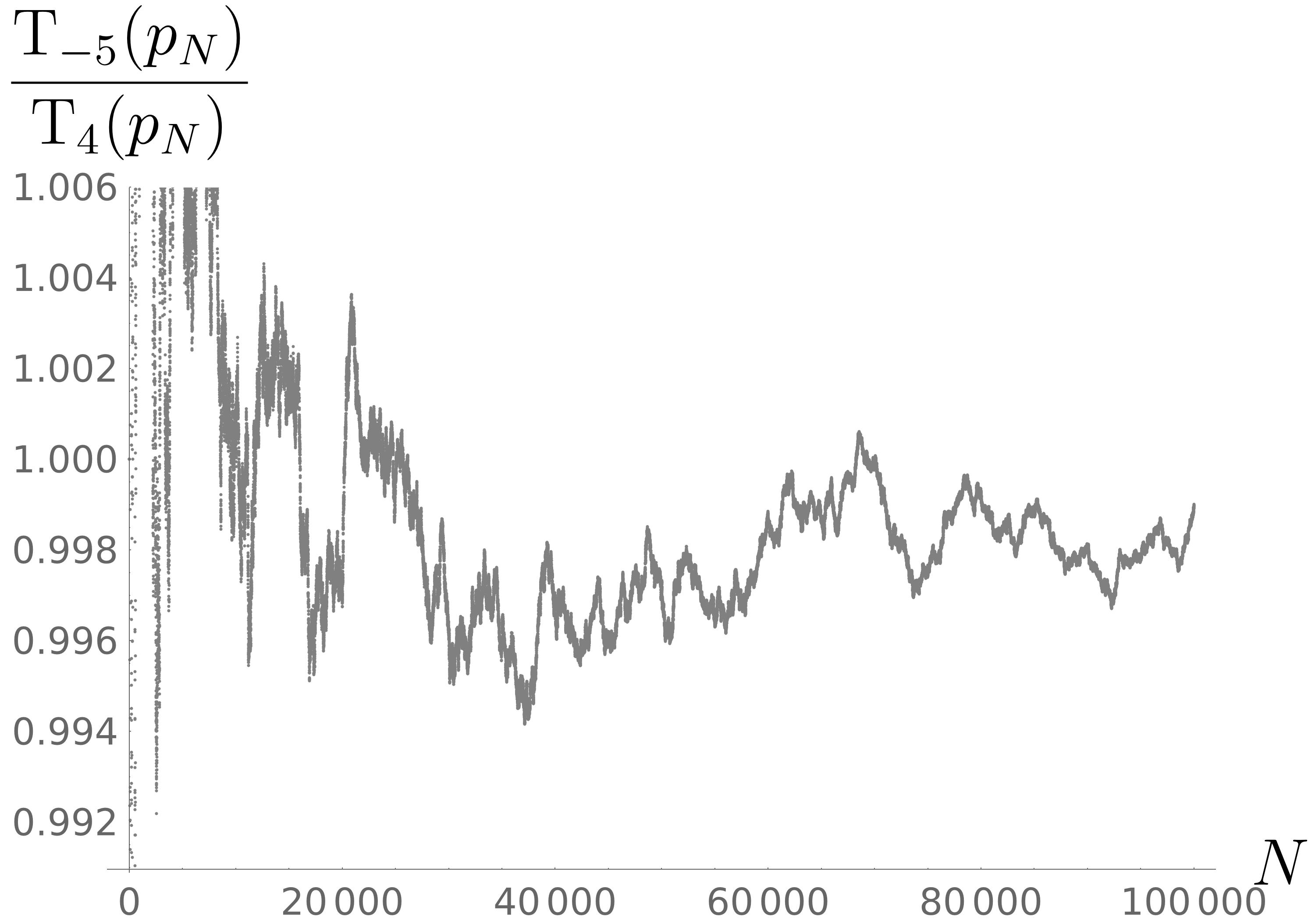}
		
		\caption{The ratio between $\ST_c(p_N)$ and $\ST_{c'}(p_N)$ for a few different $c$ and $c'$ different from $0$ and $-2$. Here $p_N$ denotes the $N$th prime number, for $N = 1, 2, \ldots, 100\,000$.}
		\label{fig:STcNequal}
	\end{figure}
	
	Using these computations of $\T_c(p)$ for consecutive primes $p$ we compute $\ST_c(N)$ for $N$ up to and including the last of these primes $p$, and observe something quite interesting: all $\ST_c(N)$, $c \neq 0, -2$ that we test appear to be asymptotically equal. We present experimental evidence of this in Figure~\ref{fig:STcNequal}, in which we look at the ratio $\ST_c(p_N) / \ST_{c'}(p_N)$ for $c \neq c'$, with $c, c' \neq 0, -2$, where $p_N$ is the $N$th prime.
	
	For the one hundred thousand data points collected for each $c = -5, -4, \ldots, 5$, it appears as though these ratios all approach 1. Indeed, though not all are shown in Figure~\ref{fig:STcNequal}, all of the plots of this type approach 1 as $N$ increases, for any combination of $c, c' \in \Set{-5, -4, -3, -1, 1, 2, 3, 4, 5}$.
	
	We observe the same behaviour, though on a smaller data set, for all combinations of $c, c' \in \Set{-100, -99, \ldots, 100} \setminus \Set{0, -2}$, for the first ten thousand primes.
	
\begin{figure}
		\centering
		\includegraphics[width = 0.48 \textwidth]{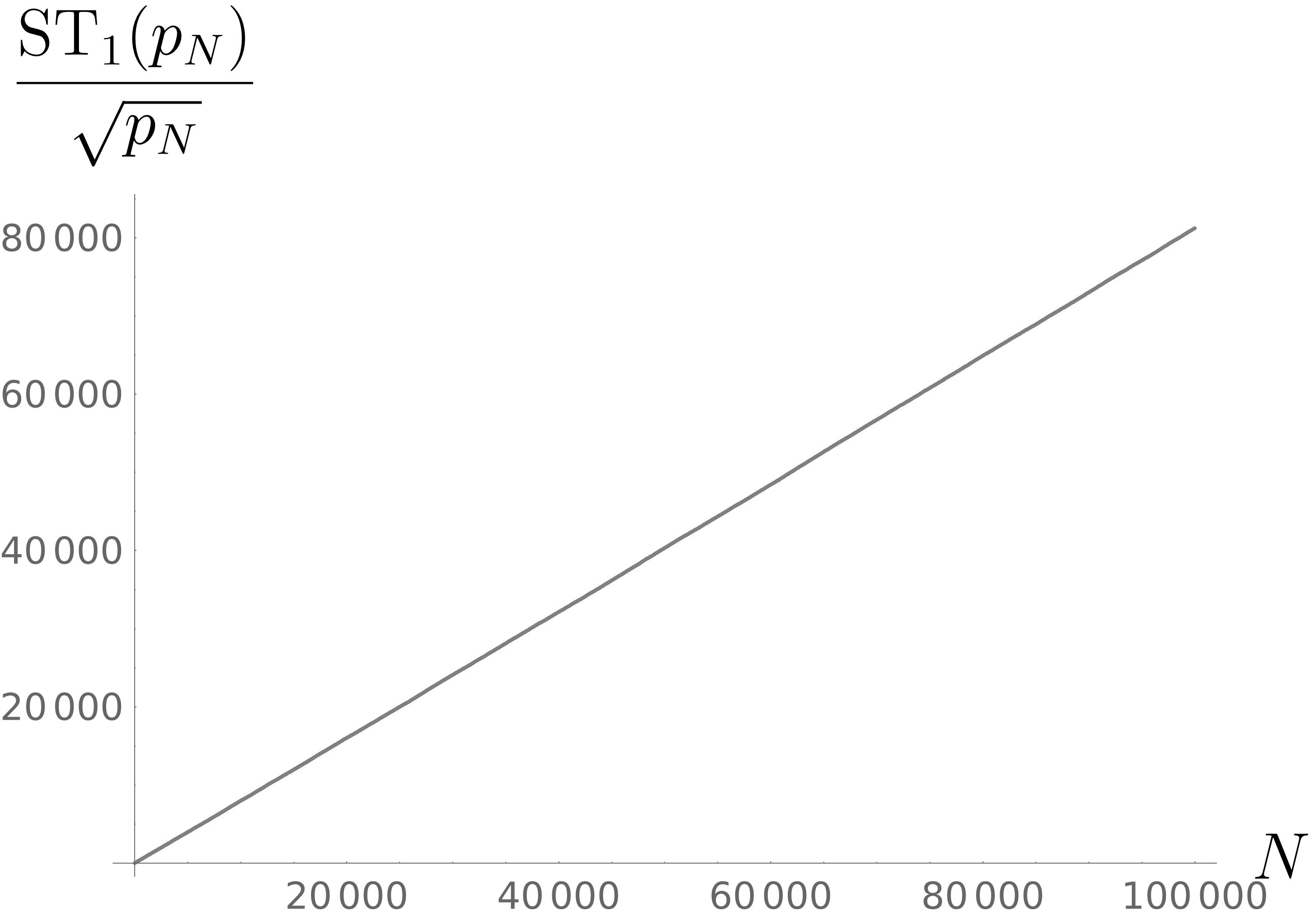} \includegraphics[width = 0.48 \textwidth]{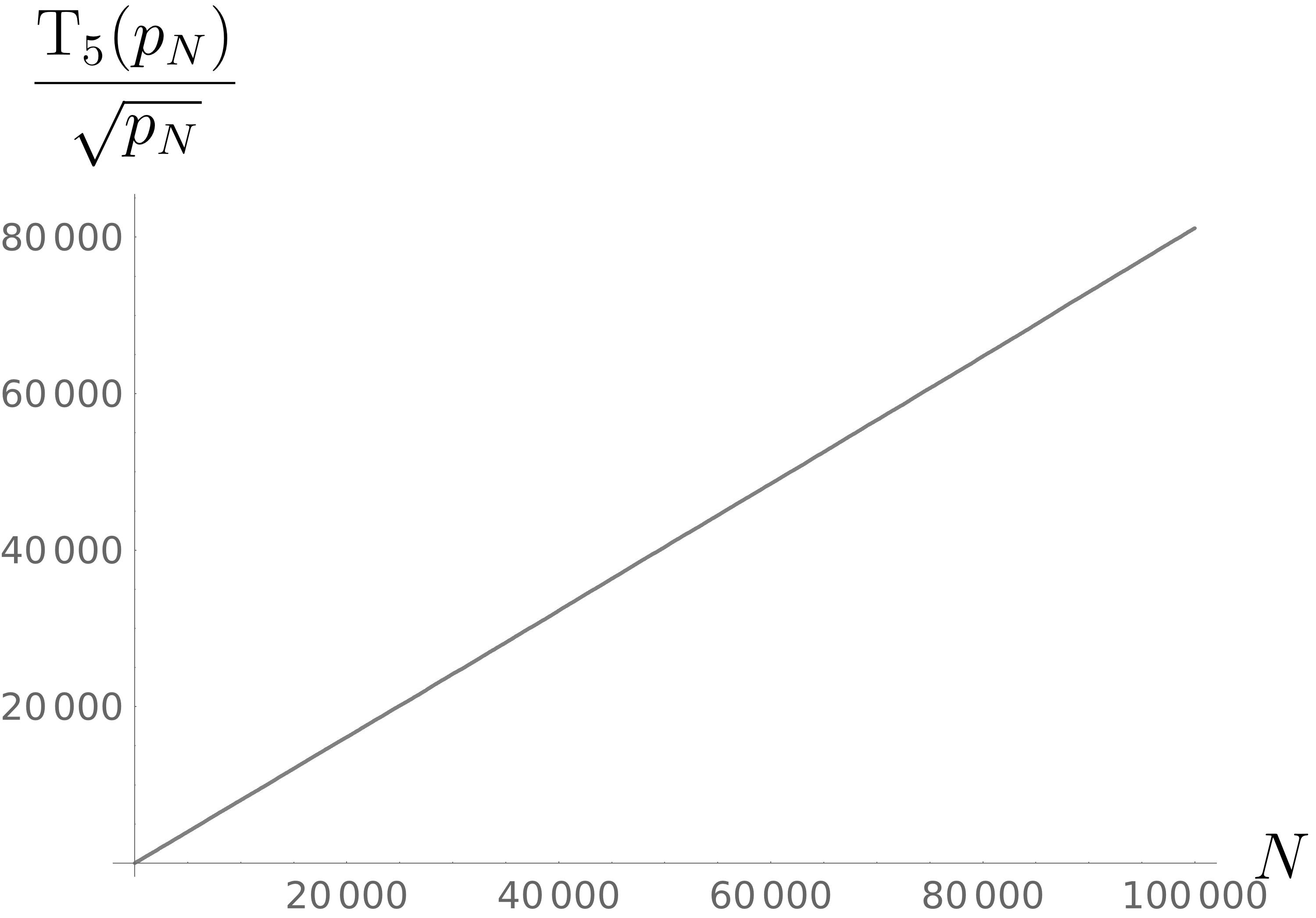}
		
		\caption{Plots of $\ST_c(p_N) / \sqrt{p_N}$ for a two different $c \neq 0, -2$, for the first $100\,000$ primes $p_N$.}
		\label{fig:STcsqrt}
	\end{figure}

	Since $\ST_c(p_N)$ and $\ST_{c'}(p_N)$ for $c, c' \neq 0, -2$ appear to be asymptotically equal, it is natural to ask what they might be equal to. We arrive at what might potentially be an answer to this question by dividing $\ST_c(p_N)$ by $\sqrt{p_N}$ and finding that this ratio appears to be linear. We show this plot for $c = 1$ and $c = 5$ in Figure~\ref{fig:STcsqrt}. Again we observe the same behaviour for all $c \in \Set{-100, -99, \ldots, 100} \setminus \Set{0, -2}$, for the first one hundred thousand primes for $c \neq 0, -2$ between $-5$ and $5$ and for the first ten thousand primes for the rest.
	
	That these plots look very similar is of course to be expected, given how previous investigations suggest that different $\ST_c(p_N)$ appear to be very close to equal for $c \neq 0, -2$ as $N$ grows large.
	
	\begin{figure}
		\centering
		\includegraphics[width = 0.48 \textwidth]{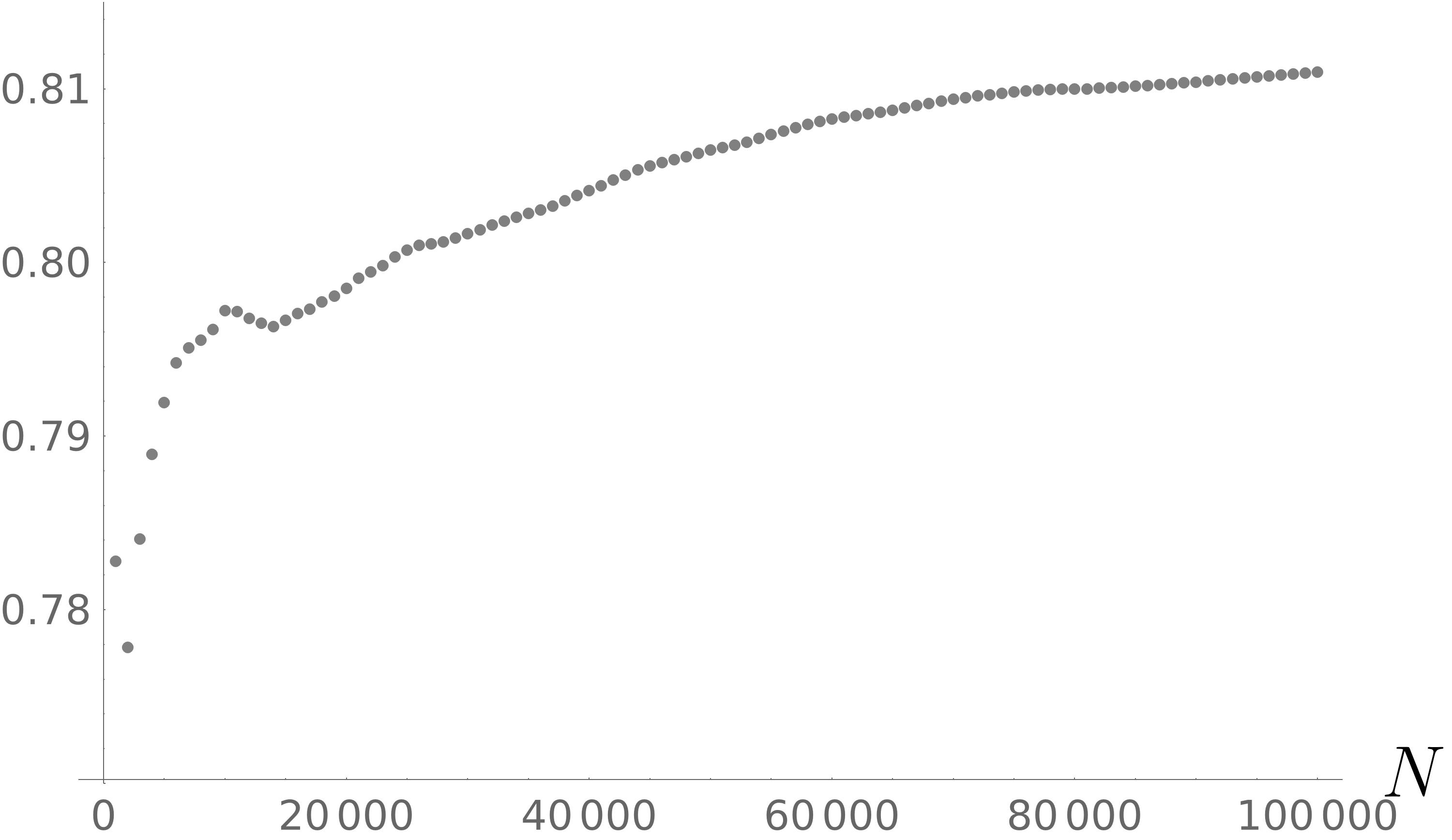} \quad \includegraphics[width = 0.48 \textwidth]{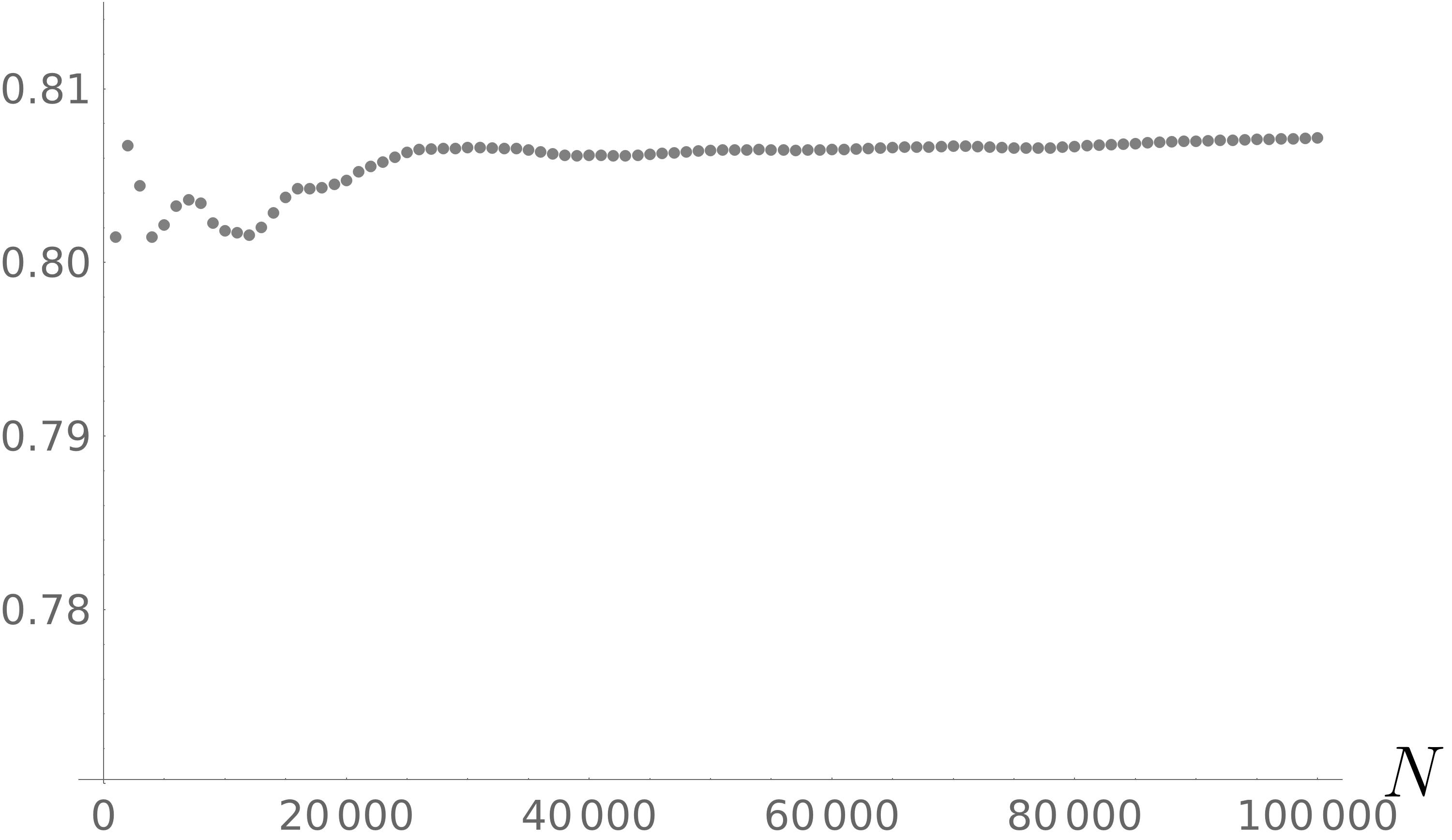}
		\includegraphics[width = 0.48 \textwidth]{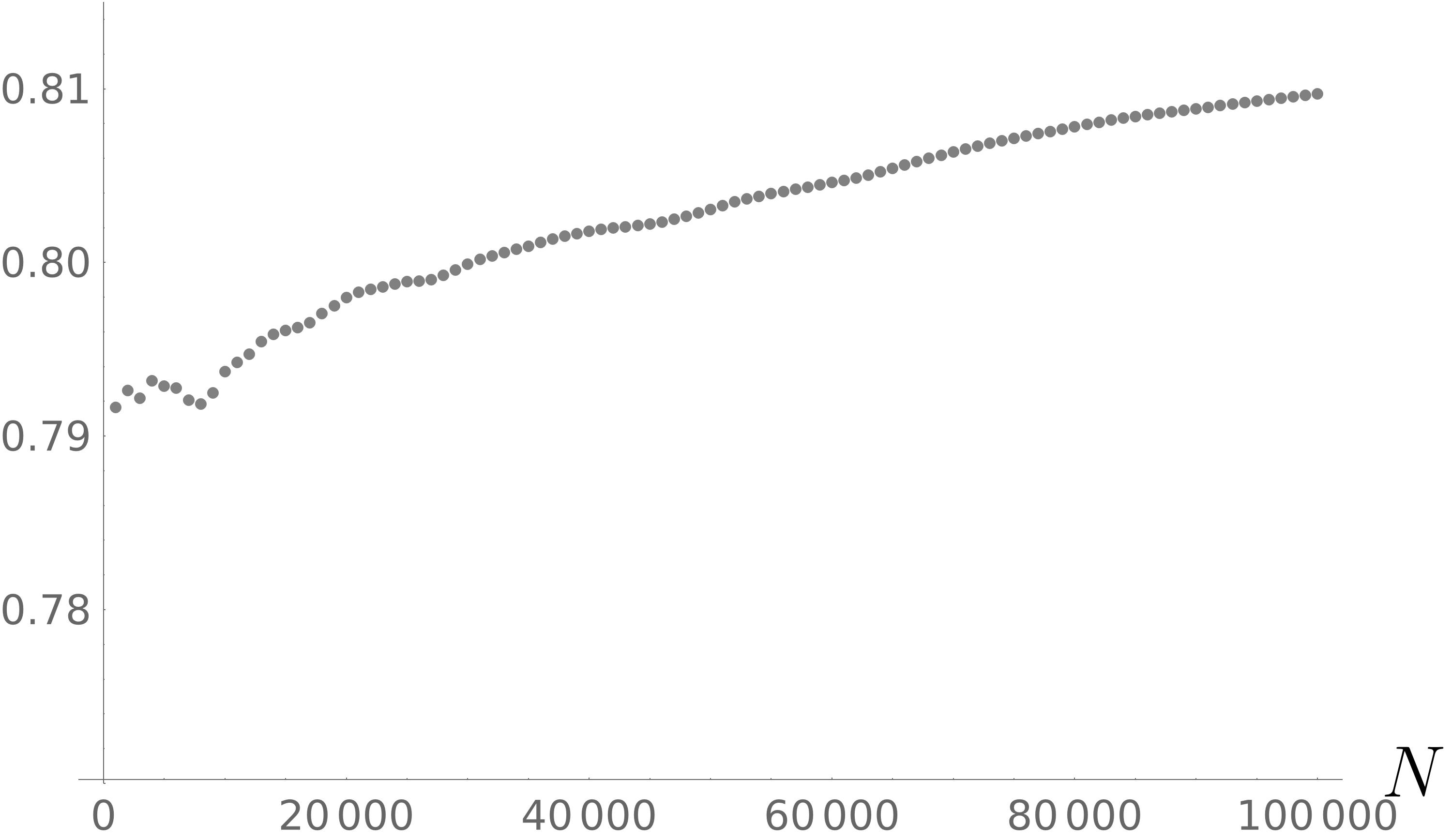} \quad \includegraphics[width = 0.48 \textwidth]{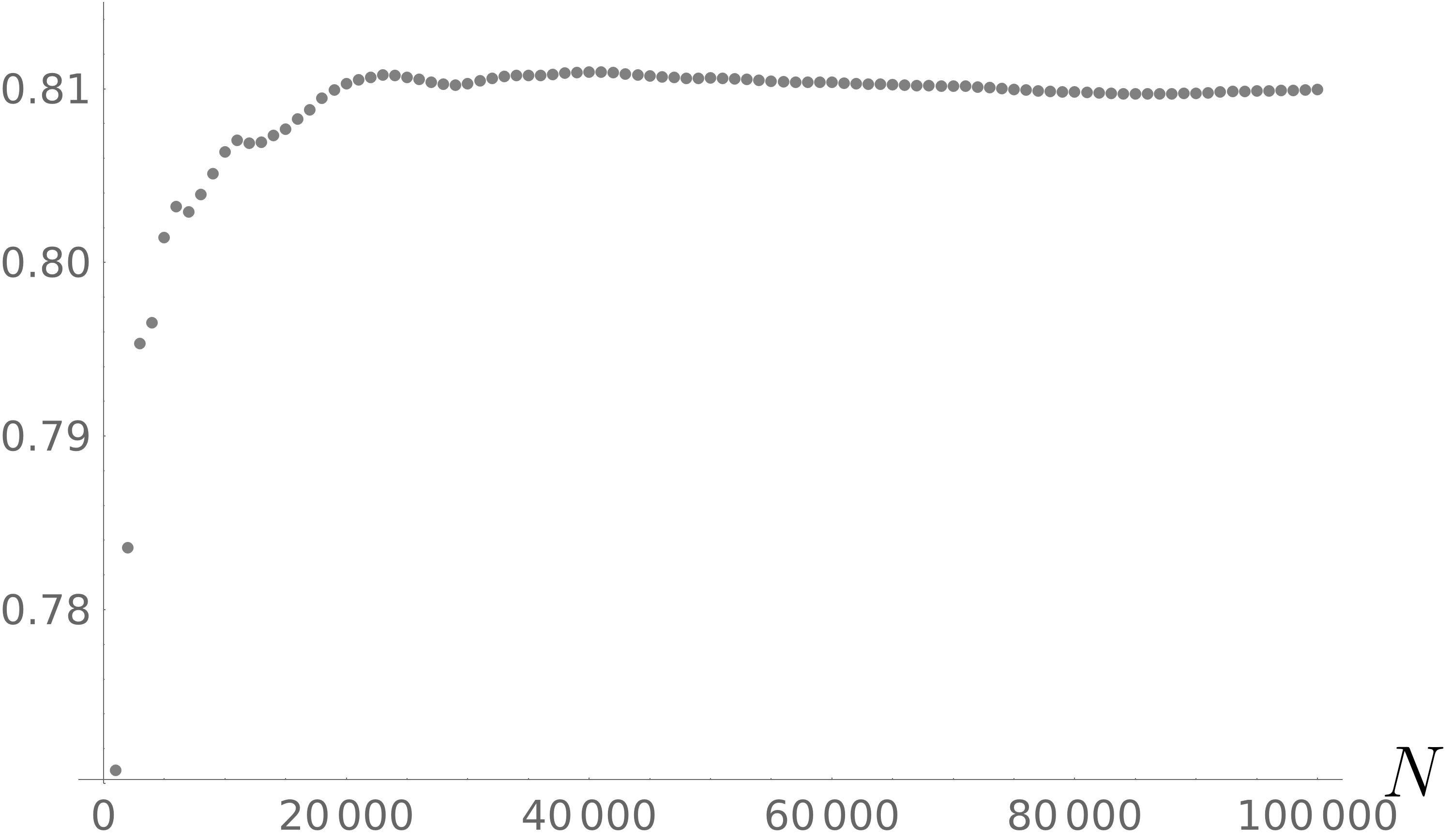}
		\caption{The slope of the linear approximations of $\ST_{c}(p_N) / \sqrt{p_N}$ for $c = -4, -5, 1, 2$ (from left to right, top to bottom) for the first $N$ primes $p_N$, for $N = 1\,000, 2\,000, 3\,000, \ldots, 100\,000$.}
		\label{fig:STcsqrtslope}
	\end{figure}
	
	Now since the plot of $\ST_c(p_N) / \sqrt{p_N}$ appears linear, we would like to decide the slope of this line. It varies slightly, although not by very much, for our different $c \neq 0, -2$ between $-5$ and $5$, the smallest being roughly $0.807$ for $c = -5$ and the largest being roughly $0.811$ for $c = -4$. They all seem to be fairly stable around $0.81$, which we observe by plotting the slope of the linear approximation for $\ST_c(p_N) / \sqrt{p_N}$ for the first $N$ primes, with $N$ taken to be $1\,000, 2\,000, 3\,000, \ldots, 100\,000$ in Figure~\ref{fig:STcsqrtslope}. 
	
	Note that Figure~\ref{fig:STcsqrt} suggets that for $c \neq 0, -2$
	\[
		\frac{\ST_c(p_N)}{\sqrt{p_N}} \sim kN,
	\]
	where $p_N$ is the $N$th prime number and $k$ is some constant that Figure~\ref{fig:STcsqrtslope} suggests is close to $0.81$. Further note that since $p_N$ is the $N$th prime number, $N$ is the number of primes less than or equal to $p_N$. We know from the Prime number theorem (see \cite{Koch1901}) that this is approximately equal to $p_n / \log p_N$, whereby we have
	\[
		\frac{\ST_c(p_N)}{\sqrt{p_N}} \sim k \frac{p_N}{\log p_N},
	\]
	which when multiplying by the square root becomes
	\[
		\ST_c(p_N) \sim k \frac{p_N^{3/2}}{\log p_N}.
	\]
	
	\noindent
	This---which it is interesting to note does not entirely unlike the asymptotic behaviour of $\ST_0(N)$ and $\ST_{-2}(N)$, shown to be $N^2 / (6 \log N)$---leads us to the formulation of our conjecture regarding the asymptotic behaviour of $\ST_c(N)$ for $c \neq 0, -2$. 
	
	\begin{conjecture}
		Let $c \neq 0, -2$ be an integer and let $f_c(x) = x^2 + c$. Further let $\T_c(p)$ denote the number of periodic points of $f_c$ over the finite field $\F_p$ of $p$ elements, $p$ being prime. Finally let 
		\[
			\ST_c(N) = \sum_{p \leq N} \T_c(p),
		\]
		the sum of $\T_c(p)$ for all primes $p$ less than or equal to $N$. Then
		\[
			\ST_c(N) \sim k \frac{N^{3/2}}{\log N},
		\]
		where $k$ is a constant approximately equal to $0.81$.
	\end{conjecture}
	
	\section{Distribution of $\T_c(p)$}
	
	Having found that $\ST_c(N)$, which is effectively averaging much of the hard to predict behaviour of $\T_c(p)$ by means of summing, we attempt to find out how $\T_c(p)$ is distributed for fixed $c \neq 0, -2$. In order to do so we first normalise by dividing by $\sqrt{p}$, and then try to fit a distribution to the probability density function we get numerically for $\T_c(p)/\sqrt{p}$. 
	
\begin{figure}
		\centering
		\includegraphics[width = 0.48 \textwidth]{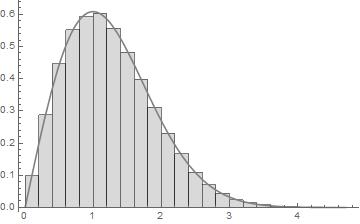} \includegraphics[width = 0.48 \textwidth]{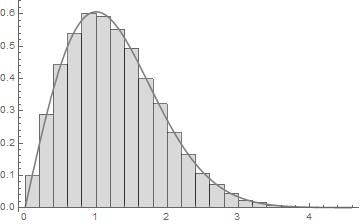}
		\caption{Plots of the histogram of $\T_c(p)/\sqrt{p}$, with $c = -5$ to the left and $c = 3$ to the right, and the probability density function of the Rayleigh distribution, with $\sigma = 0.9962$ and $\sigma = 0.9995$ respectively.}
		\label{fig:rayleigh}
	\end{figure}

	We find that for $c \neq 0, -2$ between $-100$ and $100$, $\T_c(p)$ all appear to be Rayleigh distributed with a parameter $\sigma$ between $0.978$ and $1.008$. In Figure~\ref{fig:rayleigh} we show two of the probability density histograms and the matched Rayleigh distribution. For the remaining $c$ studied they all look very similar.
	
	\section{Discussion}
	
	Given how our conjecture singles out the known cases $c = 0$ and $c = -2$ and suggests that all other $c$ behave similarly, it would be interesting to know if there are other outliers like $c = 0$ and $c = -2$. For the $c$ between $-100$ and $100$ studied, this appears to not be the case. Indeed, for these $c$, $\ST_c(p_{10\,000})$, $p_{10\,000} = 104\,729$ being the ten thousandth prime, varies between $2\,559\,211$ and $2\,637\,910$. For reference, the values of $\ST_0(p_{10\,000})$ and $\ST_{-2}(p_{10\,000})$ are both in the order of $165\,500\,000$.
	
	Further comments are warranted on the consequences of the conjecture being true. An immediate consequence is of course that $\ST_c(N)$ would behave the same for all $c \neq 0, -2$, which would be very interesting. It would mean that despite the individual $\T_c(p)$ being hard to attack there is some meaningful connection between them for different $c$. 
	
	It is also interesting that $N^2 / (6 \log N)$, the asymptotic behaviour of $\ST_0(N)$ and $\ST_{-2}(N)$, is so much greater than $k N^{3/2} / \log N$, the conjectured asymptotic behaviour for $\ST_c(N)$ for other $c$. Whilst this does not necessarily imply anything directly related to the applications discussed, such as integer factorisation and pseudorandom number generation, it does suggest that $\T_c(p)$ is usually \emph{much} smaller when $c \neq 0, -2$ than when $c = 0$ or $c = -2$. 
	
	This would mean that for a psuedorandom number generator like BBS, where we rely on large cycle lengths in order to not repeat our pseudorandom sequences too often, these maps $f_c(x) = x^2 + c$ for $c \neq 0, -2$ might not be very effective generators; $\T_c(p)$ being smaller for $c \neq 0, -2$ does not mean that all cycles are shorter, whence these generators needn't always be worse, but since we have more irregular dynamics for these $c$, we expect it to often be the case. 
	
	On another note, in \cite[Section~5]{Nilsson2013}, Nilsson looks at the number of periodic points for all $f_c$ over $\F_p$ for fixed primes $p$. That is,
	\[
		\sum_{c = 0}^{p - 1} \T_c(p).
	\]
	
	\noindent
	In doing so he plots the sum for the first one thousand primes $p$ and observes that this plot looks to behave quite well, however he does not find an explicit expression for it.
	
	Our conjecture being true might help to resolve this. In the same sense as one might use approximations of the prime counting function to estimate the number of primes in a certain interval, an approximation of $\ST_c(N)$, $c \neq 0, -2$ could be used to estimate the number of periodic points for $f_c$ over fields $\F_p$, for $p$ in some interval. Using this one can estimate the size of $\T_c(p)$ for $c \neq 0, -2$. 
	
	Using such an estimate one could approximate $p - 2$ of the terms in the above sum, together with the known explicit expressions for $\T_0(p)$ and $\T_{-2}(p)$ for the remaining two terms, and thereby describe the approximate behaviour of the sum. 
	
	We close this by mentioning that at a cursory glance it would appear as though these sums of number of periodic points for all primes less than or equal to some $N$ are well behaved not only for quadratic dynamical systems $f_c(x) = x^2 + c$, but also for larger exponents. We studied briefly the behaviour of these sums for $x \mapsto x^3 + c$, for $c$ between $0$ and $100$, where again we saw very similar behaviour for all $c \neq 0$---$c = 0$ again being the odd one out, naturally having many more periodic points---however for these maps we have not studied it to closely as to come up with a potential expression for the asymptotic behaviour of the sum.

\end{document}